\documentclass[12pt,reqno]{amsart}
\usepackage{enumerate, latexsym, amsmath, amsfonts, amssymb, amsthm, graphicx, color}
\usepackage{hyperref}
\hypersetup{hypertex=true,
	colorlinks=true,
	linkcolor=blue,
	anchorcolor=blue,
	citecolor=blue}
\usepackage{amsmath}
\def\pmod #1{\ ({\rm{mod}}\ #1)}

\def\bg{\bigg}
\def\({\bg(}
\def\){\bg)}

\def\adj{{\rm adj}}

\def\u{{\bf u}}
\def\v{{\bf v}}

\theoremstyle{plain}
\newtheorem{theorem}{Theorem}

\newtheorem{lemma}{Lemma}

\theoremstyle{definition}

\theoremstyle{remark}

\makeatletter
\@namedef{subjclassname@2010}{%
	\textup{2010} Mathematics Subject Classification}
\makeatother
\vspace{4mm}

\allowdisplaybreaks[1]
\begin{document}

	\title{On a generalization of R. Chapman's ``evil determinant''}
	
	\author[L.-Y. Wang, H.-L. Wu, H.-X. Ni]{Li-Yuan Wang, Hai-Liang Wu, He-Xia Ni*}

	\address {(Li-Yuan Wang) School of Physical and Mathematical Sciences, Nanjing Tech University, Nanjing 211816, People's Republic of China}
	\email{\tt wly@smail.nju.edu.cn}

	\address {(Hai-Liang Wu) School of Science, Nanjing University of Posts and Telecommunications, Nanjing 210023, People's Republic of China}
	\email{\tt whl.math@smail.nju.edu.cn}
	\address{(He-Xia Ni*)Department of Applied Mathematics, Nanjing Audit University, Nanjing 211815,
		People's Republic of China}
	\email{\tt nihexia@yeah.net}
	
	\begin{abstract}
  Let $p$ be an odd prime and $x$ be an indeterminate. Recently, Z.-W. Sun proposed the following conjecture:
		$$\det\left[x+\left(\frac{j-i}{p}\right)\right]_{0\le i,j\le \frac{p-1}{2}}=\begin{cases}
			(\frac{2}{p})pb_px-a_p & \mbox{if}\ p\equiv 1\pmod4,\\
			1    & \mbox{if}\ p\equiv 3\pmod4,
		\end{cases}$$
	 where $a_p$ and $b_p$ are rational numbers related to the fundamental unit and class number
	of the real quadratic field $\mathbb{Q}(\sqrt{p})$. In this paper, we confirm the above conjecture of Sun based on Vsemirnov's decomposition of  Chapman's ``evil determinant".
	\end{abstract}
	
	\thanks{
		\newline \indent
		This research was supported by the National Natural Science Foundation of China (Grant Nos. 12371004, 12201291 and 12101321) and the Natural Science Foundation of the Higher Education Institutions of Jiangsu Province (21KJB110001, 21KJB110002).
		\newline\thanks{*Corresponding author.}}
	\maketitle
	\section{Introduction}	
	\setcounter{lemma}{0}
	\setcounter{theorem}{0}
	\setcounter{corollary}{0}
	\setcounter{remark}{0}
	\setcounter{equation}{0}
	\setcounter{conjecture}{0}

	The study of determinants  involving Legendre symbols is an active topic in algebra, number theory and combinatorics. Krattenthaler summarized many results on  the evaluation of determinants  (see \cite{K1,K2}). Let $p$ be an odd prime and  $(\frac{\cdot }{p})$ be the Legendre symbol. Carlitz \cite{carlitz} showed that
	$$\det\left[\left(\frac{j-i}{p}\right)\right]_{1\le i,j\le p-1}=p^{(p-3)/2}.$$
	
In 2004, R. Chapman \cite{chapman} studied several variants of Carlitz's results. In particular, he posed seveal related conjectures. The most famous one is the ``evil determinant conjecture" (see \cite{problems,evil} for the motivation of this conjecture). When $p\equiv1\pmod4$ we let $\varepsilon_p>1$ and $h_p$ denote the fundamental unit and class number
	of the real quadratic field $\mathbb{Q}(\sqrt{p})$ respectively. And we write
\begin{align}\label{abp}
	 \varepsilon_p^{(2-(\frac{2}{p}))h_p}=a_p+b_p\sqrt{p}\
\end{align}
where $\ a_p,b_p\in\mathbb{Q}$.
 The ``evil determinant conjecture" says that
	\begin{align}\label{evil}
		 \det\left[\left(\frac{j-i}{p}\right)\right]_{0\le i,j\le \frac{p-1}{2}}=\begin{cases}
		 	-a_p & \mbox{if}\ p\equiv 1\pmod4,\\
		 	1    & \mbox{if}\ p\equiv 3\pmod4.
		 \end{cases}
	\end{align}
	 Chapman could not solve this problem for many years  so that he named it  ``evil determinant conjecture".  This conjecture was later comfirmed by Vsemirnov \cite{V3,V1} with sophisticated matrix decomposition.
	
	 In 2019, along this line, Zhi-Wei Sun \cite{S19} obtained many revelant results on this topic. For instance, if $p\equiv 1\pmod 4$ is a prime and $d\in \mathbb{Z}$, Sun showed that
	 	$$\det\left[\left(\frac{i+dj}{p}\right)\right]_{0\le i,j\le (p-1)/2}\equiv \( \(\frac{d}{p}\)d\)^{(p-1)/4} \frac{p-1}{2}! \pmod p.$$
Also, Sun \cite{S19} posed some related conjectures for further study.
Let $x$ be an indeterminate.
	Define
	$$C(x):=\det\left[x+\left(\frac{j-i}{p}\right)\right]_{0\le i,j\le \frac{p-1}{2}}.
	$$
Recently Sun posed a conjecture on $C(x)$. In this paper we obtain the following theorem, which confirmed that conjecture.
	\begin{theorem}\label{gevil}
		Let $p$ be an odd prime and $a_p$, $b_p$ be as (\ref{abp}). Then
			$$C(x)=\begin{cases}
			(\frac{2}{p})pb_px-a_p & \mbox{if}\ p\equiv 1\pmod4,\\
			1    & \mbox{if}\ p\equiv 3\pmod4.
		\end{cases}$$
	\end{theorem}
Note that when we take $x=0$ the above result coincides with (\ref{evil}).  Our idea to prove this result is mainly based on Vsemirnov's decomposition \cite{V3,V1} of $\left[\left(\frac{j-i}{p}\right)\right]_{0\le i,j\le (p-1)/2}$.
	The proof of Theorem \ref{gevil} will be given in section 2 and section 3.

	\section{Proof of Theorem \ref{gevil}: Case $p\equiv3\pmod4$ }
	\setcounter{lemma}{0}
	\setcounter{theorem}{0}
	\setcounter{corollary}{0}
	\setcounter{remark}{0}
	\setcounter{equation}{0}
	\setcounter{conjecture}{0}
	We begin with a lemma we shall use in the proof.
	\maketitle
	\begin{lemma}[The Matrix Determinant Lemma]
		Let $H$ be an $n\times n$ matrix, and let $\u,\v$ be two $n$-dimensional column vectors. Then
		\begin{equation*}
			\det (H+\u\v^T)=\det H+ \v^T \adj(H)\u,
		\end{equation*}
		where $\adj (H)$ is the adjugate  matrix of $H$ and $\v^T$ is the transpose of $\v$.	
		\label{mdl}
	\end{lemma}
	The proof of this lemma can be found in \cite{mdl}.

Let $n=(p-1)/2$ throuthout the proof. We first consider the case $p\equiv3\pmod4$.

{\noindent{\bf Proof of Theorem \ref{gevil}: Case $p\equiv3\pmod4$.}} Let $\u=(1,1,\ldots,1)^T$ be a column vector of length $n+1$ and $\v =x \cdot \u $. By Lemma \ref{mdl}
		\begin{equation*}
		\det (C(x))=\det (C)+ x\cdot \u^T \adj(C)\u.
	\end{equation*}
 In view of (\ref{evil}), it is sufficient to show that $\u^T \adj(C)\u=0$. For any matrix $A=[a_{ij}]_{0\le i,j\le m }$, we use $\adj (A)=[A_{ji}]_{0\le i,j\le m}$ to denote its adjugate matrix.
 Then  $ \adj(C)=[C_{ji}]_{0\le i,j\le \frac{p-1}{2}}$  and
  $$\u^T \adj(C)\u=\sum_{k=0}^{n}\sum_{l=0}^{n}C_{k l}=\sum_{k=0}^{(p-3)/4}\sum_{l=0}^{n}C_{k l}+C_{n-k \; n-l}.$$
If we let $M_{i,j}$ denote the $(i,j)$ minor of $C$ (i.e., the determinant of the submatrix of $C$ formed by deleting the $i$th row and $j$th column) then
$$C_{k l}=(-1)^{k+l}M_{k,l}=(-1)^{k+l}\det\left[\left(\frac{j-i}{p}\right)\right]_{0\le i,j\le n\atop   i\ne k, j\ne l},$$
and
	$$C_{n-k\; n-l}=(-1)^{n-k+n-l}M_{n-k , n-l}=(-1)^{k+l}\det\left[\left(\frac{j-i}{p}\right)\right]_{0\le i,j\le n\atop   i\ne n-k, j\ne n-l}.$$
	Note that
	\begin{align*}
	&\det\left[\left(\frac{j-i}{p}\right)\right]_{0\le i,j\le n\atop   i\ne n-k, j\ne n-l}  \\
	&=(-1)^{\lfloor \frac{p-1}{4}\rfloor }\det\left[\left(\frac{n-j-i}{p}\right)\right]_{0\le i,j\le n\atop   i\ne n-k, j\ne l} \\
	&=(-1)^{\lfloor \frac{p-1}{4}\rfloor }\cdot (-1)^{\lfloor \frac{p-1}{4}\rfloor }\det\left[\left(\frac{n-j-(n-i)}{p}\right)\right]_{0\le i,j\le n\atop   i\ne k, j\ne l} \\
	&=\det\left[-\left(\frac{j-i}{p}\right)\right]_{0\le i,j\le n\atop   i\ne k, j\ne l} \\
		&=(-1)^n \det\left[\left(\frac{j-i}{p}\right)\right]_{0\le i,j\le n\atop   i\ne k, j\ne l}\\
			&=-\det\left[\left(\frac{j-i}{p}\right)\right]_{0\le i,j\le n\atop   i\ne k, j\ne l}.
	\end{align*}
Hence
$$C_{k l}+C_{n-k\; n-l}=0$$
for any  $0\le k\le (p-3)/4$ and $0\le l\le n$. It follows that  $$\u^T \adj(C)\u=0.$$
This completes the proof of Theorem \ref{gevil} for $ p\equiv 3\pmod4$. \qed
\section{Proof of Theorem \ref{gevil}: Case $p\equiv1\pmod4$ }
\setcounter{lemma}{0}
\setcounter{theorem}{0}
\setcounter{corollary}{0}
\setcounter{remark}{0}
\setcounter{equation}{0}
\setcounter{conjecture}{0}

In this section we deal with the case $ p\equiv 1\pmod4$. We begin with the following lemmas.

\begin{lemma} {\rm(\cite[Theorem 3]{V1})}.   \label{uv}
	Let $u_i,v_j(1\le i,j\le m)$ be complex numbers such that $u_iv_j\ne -1$. Then
	\begin{align*}
		&	\det \bigg[\frac{u_i+v_j}{1+u_iv_j}\bigg]_{1\le i,j\le m}  \\
		=&\frac{1}{2}\(\prod_{i=1}^{m}(  1+u_i )(  1+v_i )+(-1)^m \prod_{i=1}^{m}(  1-u_i )(  1-v_i )\)  \\
		&\quad \times\prod_{1\le i<j\le m}^{}\(   (u_i-u_j)(v_j-v_i) \)\prod_{i=1}^{m}  \prod_{j=1}^{m}(1+u_iv_j)^{-1}.
	\end{align*}
\end{lemma}
Using essentially the same method appeared in the proof of \cite[Lemma 2]{V1}, it is easy to get the following result.
\begin{lemma}\label{sum} Let notations be as above. Then
	\begin{align*}
		\frac{1}{2} \(   \prod_{j=1}^{n}\( 1+ \(\frac{j}{p}\)\zeta^{j} \)^2+\prod_{j=1}^{n}\( 1- \(\frac{j}{p}\)\zeta^{j} \)^2 \)=(-1)^{n/2}\zeta^{n(n+1)/2} b_p p.
	\end{align*}
	
\end{lemma}

{\noindent{\bf Proof of Theorem \ref{gevil}: Case $p\equiv1\pmod4$.}} By similar methods appeared  in the proof of the case $ p\equiv 3\pmod4$, it is sufficient to show that
\begin{equation}\label{ucu}
\u^T \adj(C)\u=\(\frac{2}{p}\)pb_p.
\end{equation}
	For convenience, we introduce Vsemivnov's notations \cite{V1}.
Let $\zeta =e^{2\pi i/p}$ and define square matrices $U, V, D$ of order $n+1$ with
$(i,j)$ entry as follows:
\begin{align*}
	u_{ij}=&\frac{(\frac{i}{p})\zeta^{-j-2i}+(\frac{j}{p})\zeta^{-2j-i}}  { \zeta^{-i-j}+(\frac{i}{p})(\frac{j}{p}) } , \quad     0\le i,j\le n,\\
	v_{ij}=&\zeta^{2ij} ,\quad   0\le i,j\le n,\\
	d_{ii}=&\prod_{0\le k\le n \atop k\ne i}^{}\frac{1}{\zeta^{2i}-\zeta^{2k}}, \quad    0\le i\le n,\\
		d_{ij}=&0,  \quad   0\le  i \ne j \le n.
\end{align*}
Vsemivnov \cite[Theorem 2]{{V1}} proved that
$$C=\(\frac{2}{p}\)\sqrt{p}\zeta^{(p-1)/4}\cdot VDUDV.$$	
By this we have
\begin{align}\label{utu}
\u^T \adj(C)\u &=  \(\frac{2}{p}\)^{n}(\sqrt{p})^n\zeta^{\frac{n^2}{2}}\cdot  \u^T   \cdot 	\adj V \cdot	\adj D \cdot	\adj  U \cdot 	\adj D	\cdot \adj (V) \cdot   \u     \nonumber \\
&= (\sqrt{p})^n\zeta^{\frac{n^2}{2}}\cdot  \u^T   \cdot 	\adj V \cdot	\adj D \cdot	\adj  U \cdot 	\adj D\cdot	\adj (V) \cdot   \u. \nonumber  \\
\end{align}

For any matrix $A$, we let $A^{(k)}$ denote the matrix obtained from $A$ via replacing all entries in the $(k+1)$-th row by $1$. Then
\begin{align*}
	\u^T   \cdot 	\adj V  =&	(1,1,\ldots,1)\cdot 	\adj V \\
	=& (1,1,\ldots,1)\cdot [V_{ji}]_{0\le i,j\le n}       \\
	=&\(\sum_{k=0}^{n}V_{0k},\sum_{k=0}^{n}V_{1k},\ldots ,\sum_{k=0}^{n}V_{nk}\)     \\
	=&\(  \det V^{(0)}  ,\det V^{(1)}  ,\ldots ,\det V^{(n)}    \).
\end{align*}
Since each element of the first row of $V$ is $1$, we have  $\det V^{(0)}=\det V  $ and $ \det V^{(k)}  =0$ for $1\le k\le n$.
Thus
\begin{align*}
\u^T \cdot 	\adj V  = \det V \cdot ( 1  ,  0,0,\ldots,0  ).
\end{align*}
Recall that $D$ is a diagonal matrix with  $(i+1)$-th diagonal element $d_{ii}$. Hence
\begin{displaymath}
	\adj D=\det D \cdot D^{-1}
	=\det D\cdot 	\left( \begin{array}{cccccccc}
		\frac{1}{d_{00}} & & \  &  \\
		&	\frac{1}{d_{11}} &    &  \\
		& & \ddots & \\
		& &  &	\frac{1}{d_{nn}} \\
	\end{array} \right).
\end{displaymath}
It follows that
\begin{align*}
	\u^T \cdot 	\adj V \cdot	\adj D=	\frac{\det V \cdot \det D}{d_{00}} \cdot ( 1  ,  0,0,\ldots,0  ) .
\end{align*}
Note that $V$, $D$, $\adj V$ and $\adj D$ are symmetric matrices. Hence
$$\adj D \cdot 	\adj V \cdot \u	=\frac{\det V \cdot \det D}{d_{00}} \cdot ( 1  ,  0,0,\ldots,0  )^T.$$
Consequently,
\begin{align*}
& \u^T   \cdot 	\adj V \cdot	\adj D \cdot	\adj  U \cdot 	\adj D\cdot	\adj (V) \cdot   \u     \\
&=\frac{(\det V)^2 \cdot (\det D)^2}{d_{00}^2}\cdot ( 1  ,  0,0,\ldots,0  ) \cdot	\adj  U \cdot  ( 1  ,  0,0,\ldots,0  )^T  \\
&=\frac{(\det V)^2 \cdot (\det D)^2}{d_{00}^2}\cdot (U_{00}, U_{10},\ldots,U_{n0})   \cdot  ( 1  ,  0,0,\ldots,0  )^T  \\
&=\frac{(\det V)^2 \cdot (\det D)^2}{d_{00}^2}\cdot U_{00}.   \\
\end{align*}
By (\ref{utu}) we have
\begin{align*}
	 \u^T \adj(C)\u =(\sqrt{p})^n\zeta^{\frac{n^2}{2}}\cdot \frac{(\det V)^2 \cdot (\det D)^2}{d_{00}^2}\cdot U_{00}.
\end{align*}
Recall that
\begin{align*}
	 \frac{1}{	d_{00}^{2} }=  &\prod_{0< k\le n }^{}  (1-\zeta^{2k})^2   \\
	 	=&\prod_{k=1 }^{n}  (1-\zeta^{2k})\prod_{k=1 }^{n}  (1-\zeta^{2k})   \\
	 	=&(-1)^n \cdot \prod_{k=1 }^{n}  \zeta^{2k}  \cdot \prod_{k=1 }^{n}  (1-\zeta^{-2k}) \cdot \prod_{k=1 }^{n}  (1-\zeta^{2k})     \\
	 		=&  \zeta^{n(n+1)}  \cdot \prod_{k=1 }^{n}  (1-\zeta^{-2k}) \cdot \prod_{k=1 }^{n}  (1-\zeta^{2k})    \\
	 			=&  \zeta^{n(n+1)}  \cdot \prod_{k=1 }^{p-1}  (1-\zeta^{k})     \\
	 				=& p \zeta^{n(n+1)}.
\end{align*}
 Thus,
 \begin{align}\label{utuu00}
 	\u^T \adj(C)\u =p^{\frac{p+3}{4}} \zeta^{\frac{n^2+n}{2}}\cdot (\det V)^2 \cdot (\det D)^2\cdot U_{00}.
 \end{align}
 We next evaluate $U_{00}$. One can verify that
 \begin{align}\label{u00}
 U_{00}&=\det [u_{ij}]_{1\le i,j\le n}    \nonumber  \\
 &=\det \bigg[    \frac{(\frac{i}{p})\zeta^{-j-2i}+(\frac{j}{p})\zeta^{-2j-i}}  { \zeta^{-i-j}+(\frac{i}{p})(\frac{j}{p}) }    \bigg]_{1\le i,j\le n}     \nonumber  \\
  &=\det \bigg[    \frac{(\frac{i}{p})\zeta^{-i}+(\frac{j}{p})\zeta^{-j}}  { 1+(\frac{i}{p})(\frac{j}{p}) \zeta^{i+j}}    \bigg]_{1\le i,j\le n}     \nonumber  \\
    &=\det \bigg[    \frac{(\frac{i}{p})\zeta^{i}+(\frac{j}{p})\zeta^{j}}  { 1+(\frac{i}{p})(\frac{j}{p})\zeta^{i+j}} \cdot  \frac{1}{(\frac{i}{p})(\frac{j}{p}) \zeta^{i+j}}  \bigg]_{1\le i,j\le n}     \nonumber  \\
        &=\prod_{1\le i\le n}^{}   \frac{1}{(\frac{i}{p})  \zeta^{i}}  \cdot  \prod_{1\le j\le n}^{}   \frac{1}{(\frac{j}{p})  \zeta^{j}}\cdot \det \bigg[    \frac{(\frac{i}{p})\zeta^{i}+(\frac{j}{p})\zeta^{j}}  { 1+(\frac{i}{p})\zeta^{i}(\frac{j}{p})\zeta^{j}} \bigg]_{1\le i,j\le n}    \nonumber  \\
          & =  \zeta^{-\frac{p-1}{4}}  \cdot \det \bigg[    \frac{(\frac{i}{p})\zeta^{i}+(\frac{j}{p})\zeta^{j}}  { 1+(\frac{i}{p})\zeta^{i} (\frac{j}{p})\zeta^{j}} \bigg]_{1\le i,j\le n}  .   \\
           \nonumber
 \end{align}

Take $m=n=(p-1)/2$ and $u_i=v_i=(\frac{i}{p})\zeta^{i}$ in Lemma \ref{uv} and by Lemma \ref{sum} we deduce that
\begin{align}\label{detuv}
 &\det \bigg[    \frac{(\frac{i}{p})\zeta^{i}+(\frac{j}{p})\zeta^{j}}  { 1+(\frac{i}{p})\zeta^{i}\cdot (\frac{j}{p})\zeta^{j}} \bigg]_{1\le i,j\le n}    \nonumber  \\ 	
 		&=\frac{1}{2}\(\prod_{i=1}^{n}(  1+u_i )(  1+v_i )+ \prod_{i=1}^{n}(  1-u_i )(  1-v_i )\) \cdot\nonumber  \\
 	&\prod_{1\le i<j\le n}^{}\(   (u_i-u_j)(v_j-v_i) \)\prod_{i=1}^{n}  \prod_{j=1}^{n}(1+u_iv_j)^{-1} \nonumber  \\
 	&=(-1)^{n/2}\zeta^{n(n+1)/2} b_p p\cdot \prod_{1\le i<j\le n}^{}\(   (u_i-u_j)(u_j-u_i) \)\prod_{i=1}^{n}  \prod_{j=1}^{n}(1+u_iu_j)^{-1} \nonumber \\
 		&=(-1)^{n/2+\binom{n}{2}}\zeta^{n(n+1)/2} b_p p\cdot \prod_{1\le i<j\le n}^{}(   u_i-u_j )^2\prod_{i=1}^{n}  \prod_{j=1}^{n}(1+u_iu_j)^{-1} \nonumber \\
 			&=\zeta^{n(n+1)/2} b_p p\cdot \prod_{1\le i<j\le n}^{}(   u_i-u_j )^2\prod_{i=1}^{n}  \prod_{j=1}^{n}(1+u_iu_j)^{-1} \nonumber \\
 			&=\zeta^{n(n+1)/2} b_p p\cdot \prod_{1\le i<j\le n}^{}\(   \(\frac{j}{p}\)\zeta^{j}- \(\frac{i}{p}\)\zeta^{i} \)^2\prod_{i=1}^{n}  \prod_{j=1}^{n}\(1+ \(\frac{j}{p}\)\zeta^{j} \(\frac{i}{p}\)\zeta^{i}\)^{-1}. \nonumber \\
\end{align}
Define
\begin{equation*}
	f_1= \prod_{1\le i<j\le n}^{}\(   \(\frac{j}{p}\)\zeta^{j}- \(\frac{i}{p}\)\zeta^{i} \)
\end{equation*}
and
\begin{equation*}
	f_2= \prod_{1\le i<j\le n}^{}\( 1+  \(\frac{j}{p}\)\zeta^{j}\(\frac{i}{p}\)\zeta^{i} \).
\end{equation*}
Then
\begin{align}\label{f12}
	 &\prod_{1\le i<j\le n}^{}\(   \(\frac{j}{p}\)\zeta^{j}- \(\frac{i}{p}\)\zeta^{i} \)^2\prod_{i=1}^{n}  \prod_{j=1}^{n}\(1+ \(\frac{j}{p}\)\zeta^{j} \(\frac{i}{p}\)\zeta^{i}\)^{-1} \nonumber \\
	&=f_1^2\prod_{i=1}^{n}  \prod_{j=1}^{n}\(1+ \(\frac{j}{p}\)\zeta^{j} \(\frac{i}{p}\)\zeta^{i}\)^{-1}\nonumber \\
	&=f_1^2\prod_{1\le i=j\le n}^{}\(1+ \(\frac{j}{p}\)\zeta^{j} \(\frac{i}{p}\)\zeta^{i}\)^{-1}\prod_{1\le i<j\le n}^{}\(1+ \(\frac{j}{p}\)\zeta^{j} \(\frac{i}{p}\)\zeta^{i}\)^{-1}  \nonumber\\
	&\quad \times \prod_{1\le j<i\le n}^{}\(1+ \(\frac{j}{p}\)\zeta^{j} \(\frac{i}{p}\)\zeta^{i}\)^{-1}\nonumber \\
	&= f_1^2\  f_2^{-2} \prod_{j=1}^{n}(1+ \zeta^{2j})^{-1} .\nonumber\\
\end{align}
By \cite[Corollary 2 (4.6)]{V1} we see that
\begin{equation}\label{2j}
	\prod_{j=1}^{n}(1+\zeta^{2j})=\zeta^{n(n+1)/2} \(\frac{2}{p}\).
\end{equation}
Combining (\ref{detuv}), (\ref{f12}) with (\ref{2j}), we obtain
\begin{align}\label{detf12}
	 \det \bigg[    \frac{(\frac{i}{p})\zeta^{i}+(\frac{j}{p})\zeta^{j}}  { 1+(\frac{i}{p})\zeta^{i}(\frac{j}{p})\zeta^{j}} \bigg]_{1\le i,j\le n}  =
pb_p \(\frac{2}{p}\)\cdot     f_1^2 \cdot f_2^{-2} .  \\
\nonumber
\end{align}
By (\ref{u00}) and (\ref{detf12}) we have
\begin{align*}
U_{00}= \(\frac{2}{p}\)  pb_p \cdot  \zeta^{-\frac{p-1}{4}}\cdot    f_1^2 \cdot f_2^{-2}.
\end{align*}
Combining  (\ref{utuu00}) with the above equalities, we get
 \begin{align}\label{final1}
	\u^T \adj(C)\u =\(\frac{2}{p}\) pb_p \cdot \zeta^{\frac{(p-1)^2}{8}} \cdot p^{\frac{p+3}{4}} \cdot (\det V)^2 \cdot (\det D)^2\cdot f_1^2 \cdot f_2^{-2}.
\end{align}
By \cite[Eq. (4.10)]{V1} one may verify that

\begin{align}
	1&=   \zeta^{(n+1)(p-1)/4} \( 	\(\frac{2}{p}\) \sqrt{p} \)^{n+2} \cdot (\det V)^2 \cdot (\det D)^2  \cdot (\det G)^{-2} \cdot f_1^2 \cdot f_2^{-2}  \nonumber \\
	 &= \zeta^{-\frac{ p^2-1}{8}} \cdot p^{\frac{p+3}{4}} \cdot (\det V)^2 \cdot (\det D)^2\cdot f_1^2 \cdot f_2^{-2}\nonumber
\end{align}
where $$(\det G)^2=	\(\prod_{j=1}^{n}(\frac{j}{p})\zeta^{j} \)^2 =\zeta^{n(n+1)}=\zeta^{\frac{ p^2-1}{4}} .$$
Hence
 \begin{align}\label{final2}
\zeta^{\frac{(p-1)^2}{8}}\cdot  p^{\frac{p+3}{4}} \cdot  (\det V)^2 \cdot (\det D)^2\cdot f_1^2 \cdot f_2^{-2}
 = \zeta^{\frac{(p-1)^2}{8} +\frac{p^2-1}{8}  }=1.
\end{align}
The last equality follows from
\begin{equation}
\frac{(p-1)^2}{8} +\frac{p^2-1}{8} =p\cdot \frac{p-1}{4}	\equiv 0 \pmod p. \nonumber
\end{equation}
Now (\ref{ucu}) follows from (\ref{final1}) and (\ref{final2}).

This completes the proof of Theorem \ref{gevil}.\qed

\end{document}